\documentclass[11pt]{article}
\usepackage{amsmath}
\usepackage{epsfig}
\usepackage{verbatim}
\usepackage{amssymb}
\usepackage{amsfonts}
\usepackage{amsbsy}
\usepackage{epsfig}
\usepackage{color}

\newcommand{\al}{{\alpha}}

\newcommand{\R}{{\mathbb R}}
\newcommand{\N}{{\mathbb N}}

\newcommand{\pmax}{p_{max}}
\newcommand{\GKM}{\mathcal G_{m}\left(d_1,...,d_K\right)}
\newcommand{\gp}{\mathcal G_{m}\left(d_1,...,d_K\right)}
\newcommand{\gpo}{\mathcal G^0\left(p_1,...,p_m\right)}
\newcommand{\gpt}{\mathcal G^0\left(p_{\tau(1)},...,p_{\tau(m)}\right)}
\newtheorem{prop}{Proposition}[section]
\newtheorem{lem}[prop]{Lemma}

\newtheorem{remk}[prop]{Remark}

\newtheorem{theo}[prop]{Theorem}

\numberwithin{equation}{section}

\begin{document}

\title {On the limiting shape of Young diagrams associated with
inhomogeneous random words}
\author{Christian Houdr\'e \thanks{Georgia Institute of Technology, School of Mathematics, Atlanta, Georgia, 30332-0160, houdre@math.gatech.edu. {Research supported in part by the NSA Grant H98230-09-1-0017.}}
\and Hua Xu \thanks {Southern Polytechnic State University, Department of Mathematics,
Marietta, Georgia, 30060, hxu2@spsu.edu.} }

\maketitle

\begin{abstract}
\noindent The limiting shape of the random Young diagrams associated with an inhomogeneous random
word is identified as a multidimensional Brownian functional.  This functional is identical in law
to the spectrum of a random matrix. The Poissonized word problem is also briefly studied,
and the asymptotic behavior of the shape analyzed.

\end{abstract}



\section{Introduction}

Let $X_1, X_2,...,X_n, \dots $ be a sequence of random variables taking values in an ordered
alphabet.  The length of the longest (weakly) increasing subsequence of $X_1, X_2,...,X_n$,
denoted by $LI_n$, is the maximal $1 \le k \le n$ such that there exists an increasing sequence of integers $1\le i_1<i_2<\cdots <i_k\le n$ with $X_{i_1}\le X_{i_2}\le \cdots\le X_{i_k}$, i.e.,
$$LI_n=\max\left\{k: \exists\ 1\le i_1<i_2<\cdots<i_k\le n,\
with\ X_{i_1}\le X_{i_2}\le \cdots\le X_{i_k}\right\}.$$

When the $X_i{s}$ take their values independently and uniformly in an $m$-letter ordered alphabet,
through a careful analysis of the exponential generating function of $LI_n$, Tracy and
Widom [\ref{TrWi3}] gave the limiting distribution of $LI_n$
(properly centered and normalized) as that of the largest eigenvalue
of a matrix drawn from the $m\times m$ traceless Gaussian Unitary Ensemble (GUE).
This result, motivated by the celebrated random permutation result of Baik, Deift and
Johansson [\ref{r41}], was further extended to the
non-uniform setting  by Its, Tracy and Widom ([\ref{ItTrWi}], [\ref{ItTrWi1}]).  In that last setting, the corresponding limiting law is the maximal
eigenvalue of a direct sum of mutually independent GUEs
subject to an overall trace constraint.

A method to study the asymptotic behavior of the length of
longest increasing subsequences is through Young diagrams  ([\ref{Fu}], [\ref{St}]).
Recall that a Young diagram of size $n$ is a collection of $n$ boxes arranged
in left-justified rows,
with a weakly decreasing number of boxes from row to row.
The shape of a Young diagram is the vector $\lambda=\left(\lambda_1,\lambda_2,...,\lambda_k\right)$, where $\lambda_1\ge \lambda_2\ge ...\ge \lambda_k$ and for each $i$, $\lambda_i$ is the number of boxes in the $i$th
row while $k$ is the total number of rows of the diagram (and so $\lambda_1+\cdots+\lambda_k=n$).
Recall also that a (semi-standard) Young tableau is a Young diagram, with a filling of a positive integer
in each box, in such a way that the integers are weakly increasing along the rows and
strictly increasing down the columns. A standard Young tableau of size $n$ is a
Young tableau in which the fillings are the integers from $1$ to $n$.

Let now $[m]:=\left\{1,2,...,m\right\}$ be an $m$-letter ordered alphabet. A {\it word} of length $n$
is a mapping $W$ from $\left\{1,2,...,n\right\}$ to $\left\{1,2,...,m\right\}$, and let $[m]^n$
denotes the set of words of length $n$ with letters taken from the alphabet $\left\{1,2,...,m\right\}$.
A word is a {\it permutation} if $m=n$, and $W$ is onto. The Robinson-Schensted correspondence is a bijection between the set of words $[m]^n$ and the set of pairs of Young tableaux $\left\{(P,Q)\right\}$,
where $P$ is semi-standard with entries from  $\left\{1,2,...,m\right\}$, while $Q$
is standard with entries from $\left\{1,2,...,n\right\}$.
Moreover $P$ and $Q$ share the same shape which is a partition of $n$, and so, we do not distinguish
between shape and partition. If the word is a permutation, then $P$ is also standard.
A word $W$ in $[m]^n$ can be represented uniquely as an $m\times n$ matrix ${\bf X}_W$ with entries
\begin{equation}\label{Xword}
\left({\bf X}_W\right)_{i,j}=\mathbf 1_{W(j)=i}.
\end{equation}
\noindent The Robinson-Schensted correspondence actually gives a one to one correspondence
between the set of pairs of Young tableaux and the set of matrices whose entries are
either $0$ or $1$ and with exactly a unique $1$ in each column.
This was generalized by Knuth to the set of $m\times n$ matrices with nonnegative
integer entries. Let $\mathcal M(m,n)$ be the set of $m\times n$ matrices with nonnegative
integer entries.
Let $\mathcal P(P,Q)$ be the set of pairs of semi-standard Young tableaux $(P,Q)$ sharing
the same shape
and whose size is the sum of all the entries, where $P$ has elements
in $\{1,...,m\}$ and $Q$ has elements in $\{1,...,n\}$.
The Robinson-Schensted-Knuth (RSK) correspondence is a one to
one mapping between $\mathcal M(m,n)$ and $\mathcal P(P,Q)$.
If the matrix corresponds to a word in $[m]^n$, then $Q$ is standard.

Johansson [\ref{Joha2}], using orthogonal polynomial methods, proved that the limiting shape
of the Young diagrams, associated with homogeneous words, i.e., the iid uniform $m$-letter framework,
through the RSK correspondence, is the spectrum of the traceless $m\times m$ GUE.
Since $LI_n$ is also equal to the length of the top row of the associated Young diagrams, these
results recover those of [\ref{TrWi3}].  The permutation result is also obtained
by Johansson [\ref{Joha2}], Okounkov [\ref{Ok}] and Borodin, Okounkov and Olshanki [\ref{BoOkOl}].
More recently, for inhomogeneous words and via simple probabilistic tools, the limiting
law of $LI_n$ is given, in [\ref{HoLi}], as a Brownian functional.
Via the results of Baryshnikov [\ref{Bary}] or of Gravner, Tracy and
Widom [\ref{GrTrWi}] this functional can then be identified as a maximal eigenvalue
of a certain matrix ensemble.  For the shape of the associated
Young diagrams, the corresponding open problem is resolved below.

Let us now describe the content of the present paper. In Section 2, we list some simple
properties of a matrix ensemble, which we call generalized traceless GUE;
and relate various properties of the GUE to this generalized one.  In Section 3,
we obtain the limiting shape, of the RSK Young diagrams
associated with an inhomogeneous random word, as a multivariate Brownian functional.
In turn, this functional is identified as the spectrum of an $m\times m$ element of the
generalized traceless GUE.  Therefore, the limiting law of $LI_n$ is the largest eigenvalue
of the block of the $m \times m$ generalized traceless GUE corresponding to the most probable letters.
Finally, the corresponding Poissonized word problem is studied in Section 4.

\section{Generalized Traceless GUE}

In this section, we list, without proofs,
some elementary properties of the generalized traceless GUE.  Proofs are omitted since simple consequences
of known GUE results as exposed, for example, in [\ref{Meh}] or [\ref{AGZ}], except for the
proof of Proposition~\ref{prop2.7} which relies on simple arguments presented in the Appendix.

Recall that an element of the $m \times m$ GUE is an $m \times m$ Hermitian random matrix ${\bf G}=\left(G_{i,j}\right)_{1\le i,j\le m}$, whose entries are such that: $G_{i,i}\sim N(0,1)$, for $1\le i\le m$, $Re\left(G_{i,j}\right)\sim N(0,1/2)$ and $Im\left(G_{i,j}\right)\sim N(0,1/2)$, for $1\le i<j\le m$, and $G_{i,i}$, $Re\left(G_{i,j}\right)$, $Im\left(G_{i,j}\right)$ are mutually independent for $1\le i\le j\le m$. Now, for $m\ge 1$, $k=1,...,K$ and $d_1,...,d_K$ such that $\sum_{k=1}^Kd_k=m$, let $\GKM$ be the set of random matrices ${\bf X}$ which are direct sums of mutually independent elements of the $d_k\times d_k$
GUE, $k=1,...,K$ (i.e., ${\bf X}$ is an $m\times m$ block diagonal matrix whose $K$ blocks
are mutually independent elements of the $d_k\times d_k$ GUE, $k=1,...,K$). Let $p_1,\cdots, p_m>0$, $\sum_{j=1}^mp_j=1$, be such that the
multiplicities of the $K$ distinct probabilities $p^{(1)},...,p^{(K)}$ are respectively $d_1,...,d_K$, i.e., let $m_1=0$ and for $k=2,...,K$, let $m_k=\sum_{j=1}^{k-1}d_j$, and so $p_{m_k+1}=\cdots=p_{m_k+d_k}=p^{(k)}$, $k=1,...,K$.  The generalized $m\times m$ traceless GUE associated with the
probabilities $p_1,...,p_m$ is the set, denoted by $\gpo$, of $m\times m$ matrices ${\bf X}^0$, of the form
\begin{align}\label{f6.01}
{\bf X}^0_{i,j}=\left\{
              \begin{array}{ll}
                {\bf X}_{i,i}-\sqrt {p_{i}}\sum_{l=1}^m\sqrt{p_{l}}{\bf X}_{l,l}, & \hbox{if\ $i=j$;} \\
                {\bf X}_{i,j}, & \hbox{if\ $i\ne j$,}
              \end{array}
            \right.
\end{align}
\noindent where ${\bf X}\in \gp$.
\noindent Clearly, from (\ref{f6.01}), $\sum_{i=1}^m\sqrt{p_{i}}{\bf X}^0_{i,i}=0$.
Note also that the case $K=1$ (for which $d_1=m$) recovers the traceless GUE, whose elements
are of the form ${\bf X}-tr({\bf X}){\bf I}_m/m $, with ${\bf X}$ an element of the GUE and ${\bf I}_m$ the $m \times m$ identity matrix.

Here is an equivalent way of defining the generalized traceless GUE: let ${\bf X}^{(k)}$ be the $m\times m$ diagonal matrix such that
\begin{align}\label{f6.2}
{\bf X}^{(k)}_{i,i}=\left\{
              \begin{array}{ll}
                \sqrt {p^{(k)}}\sum_{l=1}^m\sqrt{p_{l}}{\bf X}_{l,l}, & \hbox{if\ $m_k<i\le m_k+d_k$;} \\
                0, & \hbox{otherwise,}
              \end{array}
            \right.
\end{align}

\noindent and let ${\bf X}\in \gp$. Then, ${\bf X}^0:={\bf X}-\sum_{k=1}^K{\bf X}^{(k)}\in\gpo$.

Equivalently, there is an "ensemble" description of $\gpo$.

\begin{prop}\label{GeTrEn}
${\bf X}^0\in\gpo$ if and only if ${\bf X}^0$ is distributed according to the
probability distribution
\begin{align}\label{GeneralizedtracelessEnsemble}
\mathbb P\left(d{\bf X}^0\right)=C\gamma\left(d{\bf X}_{1,1}^0,...,d{\bf X}_{m,m}^0\right)&\prod_{k=1}^K\Bigg(e^{-\underset{\tiny{\tiny m_k< i<j\le m_k+d_k}}{\sum}\left|{\bf X}_{i,j}^0\right|^2}\nonumber\\
&\underset{\tiny{m_k< i<j\le m_k+d_k}}{\prod}d\text{Re}\left({\bf X}^0_{i,j}\right)d\text{Im}\left({\bf X}^0_{i,j}\right)\Bigg),
\end{align}
\noindent on the space of $m\times m$ Hermitian matrices, which are direct sum of $d_k\times d_k$ Hermitian matrices, $k=1,...,K$, $\sum_{k=1}^Kd_k=m$, and where $m_1=0$, $m_k=\sum_{j=1}^{k-1}d_j$, $k=2,...,K$. Above, $C=\pi^{-\sum_{k=1}^Kd_k\left(d_k-1\right)/2}$ and $\gamma\left(d{\bf X}_{1,1}^0,...,d{\bf X}_{m,m}^0\right)$ is the distribution of an $m$-dimensional centered (degenerate) multivariate Gaussian law with covariance matrix
 $${\bf \Sigma}^{\ \!\!0}=\begin{pmatrix}
     {1-p_{1}} & -\sqrt{p_{1}p_{2}} & \cdots & -\sqrt{p_{1}p_{m}} \\
     -\sqrt{p_{2}p_{1}} & 1-p_{2} & \cdots & -\sqrt{p_{2}p_{m}}\\
     \vdots & \ddots & \ddots & \vdots \\
     -\sqrt{p_{m}p_{1}} & \cdots & -\sqrt{p_{m}p_{m-1}} & 1-p_{m} \\
   \end{pmatrix}.
$$
\end{prop}

We provide next a relation between the spectra of $\bf X$ and ${\bf X}^0$.

\begin{prop}\label{prop6.02}
Let ${\bf X}\in \GKM$, and let ${\bf X}^0\in\gpo$. Let $\xi_{1},\cdots,\xi_{m}$ be the eigenvalues of ${\bf X}$, where for each $k=1,...,K$, $\xi_{m_k+1},\cdots,$ $\xi_{m_k+d_k}$ are the eigenvalues of the $k$th diagonal block (an element of the  $d_k\times d_k$ GUE). Then, the eigenvalues of ${\bf X}^0$ are given by:
$$\xi_{i}^0=\xi_{i}-\sqrt {p_{i}}\sum_{l=1}^m\sqrt{p_{l}}{\bf X}_{l,l}=\xi_{i}-\sqrt {p_{i}}\sum_{l=1}^m\sqrt{p_{l}}{\xi}_{l},\ \ \ i=1,...,m.$$
\end{prop}

Let $\xi^{GUE,m}_{1},\xi^{GUE,m}_{2},...,\xi^{GUE,m}_{m}$ be the eigenvalues of an element of the $m\times m$ GUE. It is well known that the empirical distribution of the eigenvalues $\left(\xi^{GUE,m}_{i}/\sqrt m\right)_{1\le i\le m}$ converges almost surely to the semicircle law $\nu$ with density $\sqrt{4-x^2}/2\pi$, $-2\le x\le 2$. Equivalently, the semicircle law is also the almost sure limit of the empirical spectral measure for the $k$th block of the generalized traceless GUE, provided $d_k\rightarrow \infty$, $k=1,...,K$. This is, for example, the case of the uniform alphabet, where $K=1$, $d_1=m$ and $p^{(1)}=1/m$.

\begin{prop}\label{prop2.3}
Let $\xi^0_{1},\xi^0_{2},...,\xi^0_{m}$ be the eigenvalues of an element of the $m\times m$ generalized traceless GUE, such that $\xi^0_{m_k+1},\cdots,\xi^0_{m_k+d_k}$ are the eigenvalues of the $k$th diagonal block, for each $k=1,...,K$. For any $k=1,...,K$, the empirical distribution of the eigenvalues $\left(\xi_{i}^0/\sqrt {d_k}\right)_{m_k< i\le m_k+d_k}$ converges almost surely to the semicircle law $\nu$ with density $\sqrt{4-x^2}/2\pi$, $-2\le x\le 2$, whenever \ $d_k\rightarrow \infty$.
\end{prop}

Now for $p_1,...,p_m$ considered, so far, i.e., such that the multiplicities of the $K$
distinct probabilities $p^{(1)},...,p^{(K)}$ are respectively $d_1,...,d_K$ and $p_{m_k+1}=\cdots=p_{m_k+d_k}=p^{(k)}$, $k=1,...,K$, let
\begin{align}\label{f6.03}
\mathcal L^{p_1,...,p_m}:=\Bigg\{x=(x_1,...,x_m)\in\ \R^m:\ x_{m_k+1}&\ge\cdots\ge x_{m_k+d_k},\ k=1,...,K;\nonumber\\
& \sum_{j=1}^m\sqrt{p_{j}}x_j=0\Bigg\}.
\end{align}
\noindent In other words, $\mathcal L^{p_1,...,p_m}$ is a subset of the hyperplane $\sum_{j=1}^m\sqrt{p_{j}}x_j=0$, where within each block of size $d_k$, $k=1,...,K$, the coordinates $x_{m_k+1},...,x_{m_k+d_k}$, are ordered. For any $s_1,...,s_m\in \R$, let also
\begin{align}\label{f6.04}
\mathcal L^{p_1,...,p_m}_{(s_1,...,s_m)}:=\mathcal L^{p_1,...,p_m} \cap\Big\{(x_1,...,x_m)\in\ \R^m:
x_i\le s_i,\ i=1,...,m\Big\}.
\end{align}
\noindent The distribution function of the eigenvalues, written in non-increasing order
within each $d_k\times d_k$ GUE, of an element of $\gpo$ is given now.

\begin{prop}\label{prop2.4}
The joint distribution function of the eigenvalues, written in non-increasing order within each $d_k\times d_k$ GUE, of an element of $\gpo$ is given, for any $s_1,...,s_m\in \R$, by
\begin{align}\label{EigenDistribution}
\mathbb P\Big(\xi^0_1\le s_1,&\xi^0_2\le s_2,...,\xi^0_m\le s_m\Big)={\int}_{\!\!\!\!\mathcal L^{p_1,...,p_m}_{(s_1,...,s_m)}}\!\!\!\!f(x)dx_1\cdots dx_{m-1},
\end{align}
\noindent where for $x=(x_1,...,x_m)\in \R^m$,
\begin{equation}\label{f6.1}
f(x):=c_{m}\prod_{k=1}^K\Delta_{k}(x)^2e^{-\sum_{i=1}^m x^2_i/2}\mathbf 1_{\mathcal L^{p_1,...,p_m}}(x),
\end{equation}

\noindent with $c_{m}=(2\pi)^{-(m-1)/2}\prod_{k=1}^K\left(0!1!\cdots\left(d_k-1\right)!\right)^{-1}$ and where $\Delta_{k}(x)$ is the Vandermonde determinant associated with those $x_i$ for which $p_{i}=p^{(k)}$, i.e.,
$$\Delta_{k}(x)=\underset{\tiny{m_k+1\le i<j\le m_k+d_k}}{\prod}\left(x_i-x_j\right).$$

\end{prop}

\begin{remk}
 When the eigenvalues are not ordered within each $d_k\times d_k$ GUE, the identity (\ref{EigenDistribution}) remains valid, multiplying $c_m$, above, by $\prod_{k=1}^K\left(d_k!\right)^{-1}$, and also by omitting the ordering constraints $x_{m_k+1}\ge\cdots\ge x_{m_k+d_k},\ k=1,...,K$, in the definition of $\mathcal L^{p_1,...,p_m}$.
\end{remk}

The next proposition gives a relation in law between
the spectra of elements of $\GKM$ and of $\gpo$.

\begin{prop}\label{prop6.2} For any $m\ge 2$, let ${\bf X}\in\GKM$ and let ${\bf X}^0\in\gpo$. Let $\xi_{1},\cdots,\xi_{m}$ be the eigenvalues of ${\bf X}$, and let $\xi^0_{1},\cdots,\xi^0_{m}$ be the eigenvalues of
${\bf X}^0$ as given in Proposition \ref{prop6.02}.  Then,
$$\left(\xi_{1},\cdots,\xi_{m}\right)\stackrel{d}{=}\left(\xi^0_{1},\cdots,\xi^0_{m}\right)+\left(Z_1,\cdots,Z_m\right),$$
where $\left(Z_1,\cdots,Z_m\right)$ is a centered (degenerate) multivariate Gaussian vector with covariance matrix $\left(\sqrt{p_{i}p_{j}}\right)_{1\le i,j\le m}$. Moreover, $\left(\xi^0_{1},\cdots,\xi^0_{m}\right)$ and $\left(Z_1,\cdots,Z_m\right)$ can be chosen independent.
\end{prop}

The asymptotic behavior of the maximal eigenvalues, within each block, of ${\bf X}^0\in \gpo$ is well known and
well understood (see also Proposition~\ref{theo4.1} and Proposition~\ref{theo5.13} of the
Appendix for elementary arguments leading to the result below).

\begin{prop}\label{prop2.7}
For $k=1,...,K$, let $\underset{m_k< i\le m_k+d_k}{\max}\xi^0_i$ be the
largest eigenvalue of the $d_k\times d_k$ block of ${\bf X}^0\in\gpo$, then
$$\lim_{d_k\rightarrow \infty}\frac{\underset{m_k< i\le m_k+d_k}{\max}\xi^0_i}{\sqrt{d_k}}=2,$$
\noindent both almost surely and in the mean.
\end{prop}

\section{Random Young Diagrams and Inhomogeneous Words}
Throughout the rest of this paper, let $W=X_1X_2\cdots X_n$ be a random word,
where $X_1,X_2,\cdots, X_n$ are iid random variables with $\mathbb P\left(X_1= j\right)=p_j$,
where $j=1,...,m$, $p_j>0$, and $\sum_{j=1}^mp_j=1$.  Let $\tau$ be a permutation of $\left\{1,...,m\right\}$ corresponding
to a non-increasing ordering of $p_1,p_2,...,p_m$, i.e.,  $p_{\tau(1)}\ge\cdots\ge p_{\tau(m)}$.
Assume also there are $k=1,...,K$, distinct probabilities in $\left\{p_1,p_2,...,p_m\right\}$, and
reorder them as $p^{(1)}>\cdots> p^{(K)}$, in such a way that the multiplicity of each $p^{(k)}$ is $d_k$,  $k=1,...,K$.  In our notation, $K=1$ corresponds to the uniform case,
where $d_1=m$. Let $m_1=0$ and for any $k=2,...,K$, let $m_k=\sum_{j=1}^{k-1}d_j$ and
so the multiplicity of each  $p_{\tau(j)}$ is $d_k$ if $m_k<\tau(j)\le m_k+d_k$, $j=1,...,m$.
Finally, let ${\bf X}_W$ be as in (\ref{Xword}) the matrix corresponding to such a random
word $W$ of length $n$.

Its, Tracy and Widom ([\ref{ItTrWi}], [\ref{ItTrWi1}]) have obtained the limiting law of the length of the longest increasing subsequence of such a random word. To recall their result, let $\left(\xi_1,...,\xi_m\right)$ be the eigenvalues of an element of $\gpt$, written in such a way that $\left(\xi_1,...,\xi_m\right)$ $=\bigg(\xi_1^{d_1},...,\xi_{d_1}^{d_1},$ $...,\xi_1^{d_K},...,\xi_{d_K}^{d_K}\bigg)$, i.e., $\xi_1^{d_k},...,\xi_{d_k}^{d_k}$ are the eigenvalues of the $k$th block, $k=1,...,K$.
Then (see [\ref{ItTrWi1}]), the limiting law of the length of the longest increasing subsequence,
properly centered and normalized, is the law of $\underset{1\le i\le d_1}{\max}\xi_i^{d_1}$.
A representation of this limiting law, as a Brownian functional is given in [\ref{HoLi}].
A multidimensional Brownian functional representation of the whole shape of the
diagrams associated with a Markov random word is further given in [\ref{HoLi3}]
(see also Chistyakov and G\"otze [\ref{CG}] or [\ref{HoLi2}] for the binary case).
Below, we obtain the convergence of the whole shape of the diagrams, in the iid non-uniform case via a
different set of techniques which is related to the work of Glynn and Whitt [\ref{GlWh}],
Baryshnikov [\ref{Bary}],
Gravner, Tracy and Widom [\ref{GrTrWi}] and Doumerc [\ref{Do}].

Let $\big(\hat B^1(t), \hat B^2(t), ..., \hat B^m(t)\big)$ be the
$m$-dimensional Brownian motion having covariance matrix
\begin{equation}\label{CovNonuniBM}
{\bf \Sigma}_t:=\left(\begin{matrix}
p_{\tau(1)}\left(1-p_{\tau(1)}\right) & -p_{\tau(1)}p_{\tau(2)} & \cdots & -p_{\tau(1)}p_{\tau(m)} \\
     -p_{\tau(2)}p_{\tau(1)} & p_{\tau(2)}\left(1-p_{\tau(2)}\right) & \cdots & -p_{\tau(2)}p_{\tau(m)} \\       \vdots  & \vdots & \ddots  & \vdots \\ -p_{\tau(m)}p_{\tau(1)}  & -p_{\tau(m)}p_{\tau(2)} & \cdots &
p_{\tau(m)}\left(1-p_{\tau(m)}\right) \\   \end{matrix}\right)t.
\end{equation}

\noindent For each $l=1,...,m$, there is a unique $1\le k\le K$ such that $p_{\tau(l)}=p^{(k)}$, and let
\begin{equation}\label{hatLk}
\hat L_m^l=\sum_{j=1}^{m_k}\hat B^{\tau(j)}(1)+\underset{J(l-m_k,d_k)}{\sup}\sum_{j={m_k+1}}^{m_k+d_k}\sum_{i=1}^{l-m_k}\big(\hat B^{\tau(j)}(t_{j-i+1}^i)-\hat B^{\tau(j)}(t_{j-i}^i)\big),
\end{equation}
\noindent where the set $J(l-m_k,d_k)$ consists of all the subdivisions $(t_j^i)$ of $[0,1]$, $1\le i\le l-m_k$, $j\in \N$, of the form:
\begin{align}\label{f5.7555}
t_j^i\in[0,1];\ t_j^{i+1}\le t_j^i\le t_{j+1}^i&;\ t_j^i=0\ for\ j\le {m_k}\nonumber\\
 &and\ t_j^i=1\ for\ j\ge {m_{k+1}-(l-m_k)}+1.
\end{align}
\noindent With these preliminaries, we have:

\begin{theo}\label{theoNonUniform}
Let $\lambda(RSK({\bf X}_W))=\left(\lambda_1,...,\lambda_m\right)$ be the common shape of the Young diagrams associated with $W$ through the RSK correspondence. Then, as $n\rightarrow\infty$,
\begin{align}
\left(\frac{\lambda_1-np_{\tau(1)}}{\sqrt{n}},...,\frac{\lambda_m-np_{\tau(m)}}{\sqrt{n}}\right)\Longrightarrow \left(\hat L_m^1,\hat L_m^2-\hat L_m^1,...,\hat L_m^m-\hat L_m^{m-1}\right).
\end{align}
\end{theo}
\noindent {\bf Proof.} Let $\left({\bf e_j}\right)_{j=1,...,m}$ be the canonical basis of $\R^m$, and let ${\bf V}=\left(V_1,...,V_m\right)$ be the random vector such that
$$\mathbb P\left({\bf V}={\bf e_j}\right)=p_j,\ \ \ j=1,...,m.$$
\noindent Clearly, for each $1\le j\le m$,
 $$\mathbb E\left(V_{j}\right)=p_j,\ Var(V_{j})=p_j\left(1-p_j\right),$$
\noindent and for $j_1\ne j_2$, $Cov(V_{j_1},V_{j_2})=-p_{j_1}p_{j_2}$. Hence the covariance matrix of ${\bf V}$ is \begin{equation}\label{f2.1.222}
{\bf \Sigma}=\left(\begin{matrix}
p_j\left(1-p_j\right) & -p_{1}p_{2} & \cdots & -p_{1}p_{m} \\
     -p_{2}p_{1} & p_2\left(1-p_2\right) & \cdots & -p_{2}p_{m} \\       \vdots  & \vdots & \ddots  & \vdots \\ -p_{m}p_{1}  & -p_{m}p_{2} & \cdots &
p_m\left(1-p_m\right) \\   \end{matrix}\right).
\end{equation}
\noindent Let ${\bf V_1},{\bf  V_2},...,{\bf V_n}$ be independent copies of ${\bf V}$, where ${\bf V_i}=(V_{i,1},V_{i,2},...,V_{i,m})$, $i=1,...,n$. Then ${\bf X}_W$ has the same law as the matrix
formed by all the $V_{i,j}$ on the lattice $\{1,...,n\}\times \{1,...,m\}$.

It is a well known combinatorial fact (see Section 3.2 in [\ref{Fu}]) that, for all $1\le l\le m$,
\begin{align}\label{e4.2.000}
\lambda_1+\cdots+\lambda_l=G^l(m,&n):=\max\Bigg\{\sum_{(i,j)\in \pi_1\cup\cdots\cup\pi_l}V_{i,j}:\pi_1,...,\pi_l\in\mathcal P(m,n),\nonumber\\
&and\ \pi_1,...,\pi_l\ are\ all\ disjoint \Bigg\},
\end{align}

\noindent where $\mathcal P(m,n)$ is the set of all paths $\pi$ taking only unit steps up or to the right in the rectangle $\{1,...,n\}\times\{1,...,m\}$ and where, by disjoint, it is meant that any two paths do not share a common point in $\left\{1,...,n\right\}\times\left\{1,...,m\right\}$ when $V_{i,j} = 1$.
We prove next that, for any $l=1,...,m$,
\begin{equation}\label{e4.4.000}
\begin{split}
\frac{G^l(m,n)-ns_l }{\sqrt{n}}\stackrel{n\rightarrow\infty}{\Longrightarrow}\hat L_m^l,
\end{split}
\end{equation}
\noindent where $s_l=\sum_{j=1}^lp_{\tau(j)}$. For $l=1$,
\begin{equation}\label{e4.222}
G^1(m,n)=\max\left\{\sum_{(i,j)\in \pi}V_{i,j}\ ;\pi\in\mathcal P(m,n)\right\}.
\end{equation}
\noindent Moreover, each path $\pi$ is uniquely determined by the weakly increasing sequence of its $m-1$ jumps, namely $0=t_0\le t_1\le \cdots \le t_{m-1}\le 1$, such that $\pi$ is horizontal on $[\lfloor t_{j-1}n\rfloor,\lfloor t_{j}n\rfloor]\times \{j\}$ and vertical on $\{\lfloor t_{j}n\rfloor\}\times[j,j+1]$. Hence
$$G^1(m,n)=\underset{0=t_0\le t_1\le \cdots \le t_{m-1}\le t_m=1}{\sup}\sum_{j=1}^m\sum_{i=\lfloor t_{j-1}n\rfloor}^{\lfloor t_{j}n\rfloor}V_{i,j}.$$
\noindent Let $p_{max}=\max_{1\le j\le m}p_j$, $J(m)=\left\{j: p_j=\pmax\right\}\subset\left\{1,...,m\right\}$ and so $d_1=card\left(J(m)\right)$ ($J(m)$ is the set of all the most probable letters). As shown in [\ref{HoLi3}, Section 3 and 4], the distribution of $G^1(m,n)$ is very close, for large $n$, to that of a very similar expression which involves only those $V_{i,j}$ for which $j\in J(m)$. To recall this result, if
$$\hat G^1(m,n)=\underset{\tiny{\begin{matrix}
                     0=t_0\le t_1\le \cdots \le t_{m-1}\le t_m=1 \\
                     t_{j-1}=t_{j}\ for\ j\notin J(m) \\
                   \end{matrix}}
}{\sup}\sum_{j=1}^m\sum_{i=\lfloor t_{j-1}n\rfloor}^{\lfloor t_{j}n\rfloor}V_{i,j},$$
\noindent then, as $n\rightarrow \infty$,

\begin{equation}\label{hatG}
\frac{G^1(m,n)}{\sqrt n}-\frac{\hat G^1(m,n)}{\sqrt n}\stackrel{\mathbb P}{\longrightarrow} 0,
\end{equation}
\noindent i.e., as $n\rightarrow \infty$, the distribution of the maximum (over all the northeast paths) in (\ref{e4.222}) is approximately the distribution of the maximum over the northeast paths going eastbound only along the rows corresponding to the most probable letters. Now,
\begin{equation}\label{e4.333}
\frac{\hat G^1(m,n)-n\pmax }{\sqrt{n}}=\!\!\underset{\tiny{\begin{matrix}
                     0=t_0\le t_1\le \cdots\\
                     \le t_{m-1}\le t_m=1 \\
                     t_{j-1}=t_{j}\ for\ j\notin J(m) \\
                   \end{matrix}}
}{\sup}\!\!\sum_{j=1}^m\frac{\sum_{i=\lfloor t_{j-1}n\rfloor}^{\lfloor t_{j}n\rfloor}V_{i,j}-(t_j-t_{j-1})n\pmax }{\sqrt{n}}.
\end{equation}
\noindent
We next claim that, as $n\rightarrow\infty$, for any $t>0$,

$$\left(\frac{\sum_{i=1}^{\lfloor tn\rfloor}V_{i,j}-tnp_j}{\sqrt{n}}\right)_{1\le j\le m}{\Longrightarrow}\left(\tilde B^j(t)\right)_{1\le j\le m},$$

\noindent where $\left(\tilde B^j(t)\right)_{1\le j\le m}$ is an $m$-dimensional Brownian motion with covariance matrix ${\bf \Sigma}t$. Indeed, for any $t>0$, since ${\bf V_1},{\bf  V_2},...$ are independent, each with mean vector $\mathbf p=(p_1,...,p_m)$, and covariance matrix ${\bf \Sigma}$,
$$\frac{\sum_{i=1}^{\lfloor tn\rfloor}{\bf V_{i}}-tn{\mathbf p}}{\sqrt{n}}\Longrightarrow \left(\tilde B^j(t)\right)_{1\le j\le m},$$
\noindent by the central limit theorem for iid random vectors and Slutsky's lemma. Next, for any $t>s>0$, and from the independence of the ${\bf V_{i}}s$,
\begin{align}
\Bigg(&\frac{{\sum_{i=\lfloor sn\rfloor+1}^{\lfloor tn\rfloor}}{\bf V_{i}}-\lfloor(t-s)n\rfloor{\mathbf p}}{\sqrt{n}},\frac{\sum_{i=1}^{\lfloor sn\rfloor}{\bf V_{i}}-\lfloor sn\rfloor{\mathbf p}}{\sqrt{n}}\Bigg)\nonumber\\
&\Longrightarrow \left(\left(\tilde B^j(t-s)\right)_{1\le j\le m},\left(\tilde B^j(s)\right)_{1\le j\le m}\right).
\end{align}

\noindent The continuous mapping theorem allows to conclude that
\begin{align}
\Bigg(&\frac{\sum_{i=1}^{\lfloor tn\rfloor}{\bf V_{i}}-tn{\mathbf p}}{\sqrt{n}},\frac{\sum_{i=1}^{\lfloor sn\rfloor}{\bf V_{i}}-sn{\mathbf p}}{\sqrt{n}}\Bigg)\nonumber\\
&\Longrightarrow \left(\left(\tilde B^j(t)\right)_{1\le j\le m},\left(\tilde B^j(s)\right)_{1\le j\le m}\right).
\end{align}

\noindent The convergence for the time points $t_1>t_2>\cdots>t_n>0$ can be treated in a similar fashion. Thus the finite dimensional distributions converge to that of $\left(\tilde B^j(t)\right)_{1\le j\le m}$. Since tightness in $C([0,1]^m)$ is as in the proof of Donsker's invariance principle (e.g., see [\ref{Bi}]), we are just left with identifying the covariance structure of the limiting Brownian motion $\left(\tilde B^j(t)\right)_{1\le j\le m}$.  But,
\begin{align}
Cov\left(\tilde B^{j_1}(t),\tilde B^{j_2}(t)\right)&=\lim_{n\rightarrow\infty}Cov\left(\frac{\sum_{i=1}^{\lfloor tn\rfloor}V_{i,j_1}}{\sqrt{n}},\frac{\sum_{i=1}^{\lfloor tn\rfloor}V_{i,j_2}}{\sqrt{n}}\right)\nonumber\\
&=\lim_{n\rightarrow\infty}\frac{1}{n}\sum_{i=1}^{\lfloor tn\rfloor}Cov\left(V_{1,j_1},V_{1,j_2}\right)\nonumber\\
&= Cov\left(V_{1,j_1},V_{1,j_2}\right)t.
\end{align}
\noindent Hence the $m$-dimensional Brownian motion $\Big(\tilde B^j(t)\Big)_{1\le j\le m}$ has covariance matrix ${\bf \Sigma}t$ with ${\bf \Sigma}$ given in (\ref{f2.1.222}). In particular, as $n\rightarrow\infty$, for any $t>0$,

$$\left(\frac{\sum_{i=1}^{\lfloor tn\rfloor}V_{i,j}-tn\pmax }{\sqrt{n}}\right)_{1\le j\le m,\ j\in J(m)}{\Longrightarrow}\left(\hat B^j(t)\right)_{1\le j\le m,\ j\in J(m)}.$$
\noindent It is also straightforward to see that the covariance matrix
of $\left(\hat B^j(t)\right)_{j\in J(m)}$ is the $d_1\times d_1$ matrix
\begin{equation}
\left(\begin{matrix}
\pmax\left(1-\pmax\right) & -\pmax^2 & \cdots & -\pmax^2 \\
     -\pmax^2 & \pmax\left(1-\pmax\right) & \cdots & -\pmax^2 \\       \vdots  & \vdots & \ddots  & \vdots \\ -\pmax^2  & -\pmax^2 & \cdots &
\pmax\left(1-\pmax\right) \\   \end{matrix}\right)t.
\end{equation}
\noindent By the continuous mapping theorem,
\begin{equation}\label{K1}
\begin{split}
\frac{\hat G^1(m,n)-np_{max} }{\sqrt{n}}\stackrel{n\rightarrow\infty}{\Longrightarrow}\underset{J(1,d_1)}{\sup}\sum_{j={1}}^{d_1}\big(\hat B^{\tau(j)}(t_{j})-\hat B^{\tau(j)}(t_{j-1})\big),
\end{split}
\end{equation}
\noindent and the right hand side of (\ref{K1}) is exactly $\hat L_m^1$, then (\ref{hatG}), leads to
\begin{equation}\label{e4.444}
\frac{G^1(m,n)-n\pmax }{\sqrt{n}}\stackrel{n\rightarrow \infty}{\Longrightarrow}\hat L_m^1.
\end{equation}

Now, for $l\ge 2$, $G^l(m,n)$ is the maximum, of the sums of the $V_{i,j}$, over $l$ disjoint paths. Still by the argument in [\ref{HoLi3}], $\left(\!G^l(m,n)\!-\!\hat G^l(m,n)\!\right)\!/\!\sqrt n\stackrel{\mathbb P}{\longrightarrow}0$, as $n\rightarrow \infty$, where $\hat G^l(m,n)$ is the maximal sums of the $V_{i,j}$ over $l$ disjoint paths we now describe. Let $1\le k\le K$ be the unique integer such that $p_{\tau(l)}=p^{(k)}$. Denote by $\al_{j(1)},...,\al_{j(m_k)}$ the letters corresponding to the $m_k$ probabilities that are strictly larger than $p_{\tau(l)}$. For each $1\le s\le m_k$, the horizontal path from $(1,j(s))$ to $(n,j(s))$ is included, and thus so are these $m_k$ paths. The remaining $l-m_k$ disjoint paths only go eastbound along the rows corresponding to the $d_k$ letters having probability $p_{\tau(l)}$. The set of these $l-m_k$ paths is in a one to one correspondence with the set of subdivisions of $[0,1]$ given in (\ref{f5.7555}). Therefore
\begin{align}
&\hat G^l(m,n)=\sum_{j=1}^{m_k}\sum_{i=1}^nV_{i,\tau(j)}+ \underset{J(l-m_k,d_k)}{\sup}\sum_{j={m_k+1}}^{m_k+d_k}\sum_{i=1}^{{l-m_k}}\sum_{r=\lfloor t_{j-i}^i  n \rfloor}^{\lfloor t_{j-i+1}^i  n\rfloor}V_{r,{\tau(j)}}.
\end{align}

\noindent Now,
\begin{align}\label{5.87}
&\frac{\hat G^l(m,n)-ns_l}{\sqrt{n}}=\sum_{j=1}^{m_k}\frac{\sum_{i=1}^nV_{i,\tau(j)}-np_{\tau(j)}}{\sqrt n}\nonumber\\
&\ \ \ \ \ \ \ \ \ \ \ \  + \underset{J(l-m_k,d_k)}{\sup}\sum_{j={m_k+1}}^{m_k+d_k}\sum_{i=1}^{l-m_k}\frac{{\sum_{r=\lfloor t_{j-i}^i n \rfloor}^{\lfloor t_{j-i+1}^i  n\rfloor}}V_{r,\tau(j)}-\left(t_{j-i+1}^i-t_{j-i}^i\right)np^{(k)}}{\sqrt n}.
\end{align}

\noindent Since the column vectors ${\bf V}_{\bf 1},{\bf  V}_{\bf 2},...,{\bf V}_{\bf n}$ are iid, again, as $n\rightarrow\infty$, for any $t>0$,

$$\left(\frac{\sum_{r=1}^{\lfloor tn\rfloor}V_{r,\tau(j)}-tnp_{\tau(j)}}{\sqrt{n}}\right)_{1\le j\le m}{\Longrightarrow}\left(\hat B^j(t)\right)_{1\le j\le m},$$

\noindent where $\left(\hat B^j(t)\right)_{1\le j\le m}$ is an $m$-dimensional Brownian motion with covariance matrix given in (\ref{CovNonuniBM}). Hence, (\ref{5.87}) and standard arguments give
$$\frac{G^l(m,n)-ns_l }{\sqrt{n}}\stackrel{n\rightarrow\infty}{\Longrightarrow}\hat L_m^l.$$

\noindent Finally, by the Cram\'er-Wold theorem, as $n\rightarrow\infty$,
\begin{align}\label{f5.82222}
\Bigg(\frac{\lambda_1-ns_1 }{\sqrt{n}},\frac{\sum_{j=1}^2\lambda_j-ns_2 }{\sqrt{n}},...,\frac{\sum_{j=1}^m\lambda_j-ns_m }{\sqrt{n}}\Bigg)\Longrightarrow \left(\hat L_m^1,\hat L_m^2,...,\hat L_m^m\right),
\end{align}

\noindent therefore, as $n\rightarrow\infty$, by the continuous mapping theorem,
\begin{align}\label{f5.82333}
&\Bigg(\frac{\lambda_1-np_{\tau(1)}}{\sqrt{n}},\frac{\lambda_2-np_{\tau(2)}}{\sqrt{n}},...,
\frac{\lambda_m-np_{\tau(m)}}{\sqrt{n}}\Bigg)\nonumber\\
&=\!\Bigg(\!\frac{G^1\!-\!ns_1 }{\sqrt{n}},\frac{\left(G^2\!-\!ns_2 \right)\!-\!\left(G^1\!-\!ns_1 \right)}{\sqrt{n}},...,\frac{\left(G^m\!-\!ns_m \right)\!-\!\left(G^{m-1}\!-\!ns_{m-1} \right)}{\sqrt{n}}\!\Bigg)\nonumber\\
&\Longrightarrow \left(\hat L_m^1,\hat L_m^2-\hat L_m^1,...,\hat L_m^m-\hat L_m^{m-1}\right).
\end{align}
\noindent The proof is now complete. \hfill $\Box$

\begin{remk}
\begin{itemize}
\item[(i)] In Theorem 3.2 of [\ref{HoLi3}], the limiting shape of the Young diagrams generated by 
an irreducible, aperiodic, homogeneous Markov word with finite state space is obtained 
as a multivariate Brownian functional similar to the one obtained above.  
The arguments there are based on a careful analysis of the reconfiguration of disjoint subsequences.  
Specifically, the smallest letter appearing in the disjoint subsequences 
is then solely in the first subsequence, the second smallest letter, not included in the 
first subsequence, is completely in the second subsequence, etc.   With this new configuration of 
the disjoint subsequences, a subdivision of the interval $[0,1]$ can be described and a Brownian functional representation is then available.  Our approach takes advantage of the lattice with zeros and ones entries 
(exactly a unique one in each column), and the fact that each subsequence corresponds to a north-east path on the lattice, and that the length of the subsequence is identical to the sum of all the 
entries on that path.  Moreover, for $1\le l\le m$, and $1\le i\le l$, the $i$th lowest path 
can be chosen to be from $(1,i)$ to $(N,M-l+i)$. Then the subdivision of $[0,1]$ is naturally determined by describing the jumps of all the paths involved.

\item[(ii)] Let $\left(\xi^0_1,...,\xi^0_m\right)$ represent the vector of the eigenvalues of an element of $\gpt$, written in such a way that $\xi^0_{m_k+1}\ge \cdots\ge \xi^0_{m_k+d_k}$ for $k=1,...,K$. Its, Tracy and Widom [\ref{ItTrWi}] have shown that the limiting density of $\left(\left(\lambda_1\!-\!np_{\tau(1)}\right)\!/\!{\sqrt{np_{\tau(1)}}},
    ...,\left(\lambda_m\!-\!np_{\tau(m)}\right)\!/\!{\sqrt{np_{\tau(m)}}}\right)$, as $n\rightarrow \infty$, is the joint density, of the eigenvalues of an element of $\gpt$, given by (\ref{f6.1}). By a simple Riemann integral approximation argument, it follows that

$$\left(\frac{\lambda_1-np_{\tau(1)}}{\sqrt{np_{\tau(1)}}},...,\frac{\lambda_m-np_{\tau(m)}}{\sqrt{np_{\tau(m)}}}\right)\Longrightarrow \left(\xi^0_1,...,\xi^0_m\right).$$
\noindent Thus, from Theorem \ref{theoNonUniform},
\begin{equation}\label{NonUniform}
\left(\frac{\hat L_m^1}{\sqrt{p_{\tau(1)}}},\frac{\hat L_m^2-\hat L_m^1}{\sqrt{p_{\tau(2)}}},...,\frac{\hat L_m^m-\hat L_m^{m-1}}{\sqrt{p_{\tau(m)}}}\right)\stackrel{d}{=}\left(\xi^0_1,...,\xi^0_m\right).
\end{equation}

\item[(iii)] Let $\big(B^1(t), B^2(t), ..., B^m(t)\big)$ be a standard $m$-dimensional Brownian motion. For $k=1,...,m$, let
$$D_m^k=\sup\sum_{i=1}^{m}\sum_{p=1}^{k}\big(B^i(t_{i-p+1}^p)-B^i(t_{i-p}^p)\big),$$
\noindent where the $\sup$ is taken over all the subdivisions $(t_i^p)$ of $[0,1]$ described in (\ref{f5.7555}).
\noindent The very approach to prove Theorem \ref{theoNonUniform} can be used to obtain a Brownian functional representation of the spectrum of the $m\times m$ GUE, namely,
\begin{equation}\label{5.34}
\left(D_m^1,D_m^2-D_m^1,...,D_m^m-D_m^{m-1}\right)
\stackrel{d}{=}\left(\xi^{GUE,m}_{1},\xi^{GUE,m}_{2},...,\xi^{GUE,m}_{m}\right).
\end{equation}
\noindent From the observation that the supremum in the definition of $G^k(m,n)$ is attained on a particular set of $k$ disjoint northeast paths for each $k=1,...,m$, Doumerc ([\ref{Do}]) found Brownian functional representations for $\sum_{i=1}^k\xi^{GUE,m}_{i}$. These functionals are similar to the $D_m^k$ except that the supremum is taken over a different set of subdivisions of $[0,1]$. In fact, we believe that the subdivisions given in (\ref{f5.7555}) should be the ones present in [\ref{Do}] (we believe the conditions $t_1 \le s_2, t_2 \le s_3, \dots$, present at the top of page 7 of [\ref{Do}], should not be there). With a similar consideration of $k$ disjoint increasing subsequences, a specific expression for the sum of the first $k$ rows of the Young diagram associated with a Markov random word is obtained, in [\ref{HoLi3}], in terms of the number of occurrences of the letters
among the sequence (see also Chistyakov and G\"otze [\ref{CG}] or [\ref{HoLi2}] for the binary case). The multidimensional convergence of the whole diagram towards a corresponding multidimensional Brownian functional is also obtained there.

In contrast to the approach in [\ref{Do}], our potential proof of (\ref{5.34}) does not require passing through the matrix central limit theorem. To briefly describe the approach in [\ref{Do}], let the $V_{i,j}$ in (\ref{e4.2.000}) be iid geometric random variables, i.e., for $r=0,1,...$, let $\mathbb P\left(V_{i,j}=r\right)=q(1-q)^r$. With such $\left\{V_{i,j}\right\}$, the probability of a given matrix realization only depend on the sum of the matrix entries, which is also the sum of the entries in the shape of the associate Young diagrams. The joint probability mass function of the shape of the associate Young diagrams through the RSK correspondence can then be expressed through the well known number of Young diagrams sharing this given shape. Next, by setting $q=1-L^{-1}$, and letting $L\rightarrow\infty$, the random variables on the lattice converge to iid exponential random variables with parameter one, while the corresponding shape of the associated Young diagrams converges to the spectrum of the $m\times n$ Laguerre Unitary Ensemble. As $n\rightarrow\infty$, for any $k=1,...,m$, the corresponding $G^k(m,n)$, properly normalized, converge in distribution to $D_m^k$. With the same normalization, it is proved in [\ref{Do}] that the spectrum of the $m\times n$ Laguerre Unitary Ensemble converges to the spectrum of the $m\times m$ GUE. Hence, the continuous mapping theorem, gives $\sum_{j=1}^k\xi^{GUE,m}_{j}\stackrel{d}{=}D_m^k$. Via the large $n$ asymptotics of the corresponding numbers of Young diagrams, we are able to directly show that the limiting joint probability mass function of the shape of the diagrams converges to the joint probability density function of the eigenvalues of an element of the GUE. Thus, $\sum_{j=1}^k\xi^{GUE,m}_{j}\stackrel{d}{=}D_m^k$, and (\ref{5.34}) follows from the Cram\'er-Wold theorem. Similar ideas are already developed by Johansson (Theorem 1.1 in [\ref{Joha2}]) to prove that the Poissonized Plancherel measure can be obtained as a limit of the Meixner measure. Johansson also proves the convergence of the whole diagram corresponding to a random word for uniform alphabets, and obtains the joint density of the limiting law.
\end{itemize}
\end{remk}

\section{The Poissonized Word Problem}

"Poissonization" is another useful tool in dealing with length asymptotics for longest
increasing subsequence problems. It was introduced by Hammersley in [\ref{Ha}] in order to show the existence of $\lim_{n\rightarrow\infty}\mathbb E L\sigma_n/\sqrt n$, for $\sigma_n$ a random permutation of $\left\{1,2,...,n\right\}$.  Since then, this technique has been widely used and we use it
below in connection with the inhomogeneous word problem.

Johansson [\ref{Joha2}] studied the Poissonized measure on the set of shapes of Young diagrams
associated with the homogeneous random word, while, Its Tracy and Widom [\ref{ItTrWi1}] also
studied the Poissonization of $LI_n$ for inhomogeneous random words.  They showed that the
Poissonized distribution of the length of the longest increasing subsequence, as a
function of $p_1,...,p_m$, can be identified as the solution of a certain integrable system of nonlinear PDEs.  Below, we show that the Poissonized distribution of the shape of the whole Young diagrams associated with an inhomogeneous random word converges to the spectrum of the corresponding direct sum of GUEs. Next, using this result, together with "de-Poissonization", we obtain the asymptotic behavior of the shape of the diagrams.

Let $W=X_1X_2\cdots X_n$ be a random word of length $n$, with each letter independently drawn and with
$\mathbb P_m\left(X_i=j\right)=p_j$, $i=1,...,n$, where $p_j>0$ and $\sum_{j=1}^mp_j=1$, i.e., the random word is distributed according to $\mathbb P_{W,m,n}=\mathbb P_m\times\cdots\times\mathbb P_m$ on the set of words $[m]^n$. Using the terminology of [\ref{Joha2}], with $\N=\left\{0,1,2,\cdots\right\}$, let
$$\mathcal P_m^{(n)}:=\left\{\lambda=\left(\lambda_1,...,\lambda_m\right)\in\N^m:\lambda_1\ge\cdots\ge\lambda_m,\ \sum_{i=1}^m\lambda_i=n\right\},$$
\noindent denote the set of partitions of $n$, of length at most $m$. The RSK correspondence defines a bijection from $[m]^n$ to the set of pairs of Young diagrams $(P,Q)$ of common shape $\lambda\in\mathcal P_m^{(n)}$, where $P$ is semi-standard with elements in $\left\{1,...,m\right\}$ and $Q$ is standard with elements in $\left\{1,...,n\right\}$.

For any $W\in [m]^n$, let $S(W)$ be the common shape of the Young diagrams associated with $W$ by the RSK correspondence. Then $S$ is a mapping from $[m]^n$ to $\mathcal P_m^{(n)}$, which, moreover, is a surjection.
The image (or push-forward) of $\mathbb P_{W,m,n}$ by $S$ is the measure $\mathbb P_{m,n}$ given, for any $\lambda_0\in\mathcal P_m^{(n)}$, by
$$\mathbb P_{m,n}\left(\lambda_0\right):=\mathbb P_{W,m,n}\left(\lambda\left(RSK({\bf X}_W)\right)=\lambda_0\right).$$
\noindent Next, let
$$\mathcal P_m:=\left\{\lambda=\left(\lambda_1,...,\lambda_m\right)\in \N^{m}:\lambda_1\ge\cdots\ge\lambda_m\right\},$$
\noindent be the set of partitions, of elements of $\N$, of length at most $m$.
The set $\mathcal P_m$ consists of the shapes of the Young diagrams associated with
the random words of any finite length made up from the $m$ letter alphabet.

For $\al>0$, the Poissonized measure of $\mathbb P_{m,n}$ on the set $\mathcal P_m$ is then defined as
\begin{equation}\label{Poissonized}
\mathbb P_m^\al\left(\lambda_0\right):=e^{-\al}\sum_{n=0}^\infty\mathbb P_{m,n}\left(\lambda_0\right)\frac{\al^n}{n!}.
\end{equation}

\noindent The Poissonized measure $\mathbb P_m^\al$ coincides with the distribution of the shape of the Young diagrams associated with a random word whose length is a Poisson random variable with mean $\al$. Such a random word is called {\it Poissonized}, and $LI_\al$ denote the length of its longest increasing subsequence.

The Charlier ensemble is closely related to the Poissonized word problem. It is used
by Johansson [\ref{Joha2}] to investigate the asymptotics of $LI_n$ for finite uniform alphabets.
For the non-uniform alphabets we consider, let us define the generalized Charlier ensemble to be:
\begin{equation}\label{GenCharlier}
\mathbb P_{Ch,m}^\al\left(\lambda^0\right)=\underset{1\le i<j\le m}{\prod}(\lambda^0_i-\lambda^0_j+j-i)\prod_{j=1}^{m}\frac{1}{(\lambda^0_j+m-j)!}s_{\lambda^0}(p)e^{-\al}\prod_{i=1}^m \al^{\lambda^0_i},
\end{equation}
\noindent for all $\lambda^0=(\lambda^0_1,\lambda^0_2,...,\lambda^0_m)\in\mathcal P_m$,
and where $s_{\lambda^0}(p)$ is the Schur function of shape $\lambda^0$ in the variable $p=\left(p_{\tau(1)},...,p_{\tau(m)}\right)$ which we describe next.
Let $\mathcal A_1,...,\mathcal A_K$ be the decomposition of $\left\{1,...,m\right\}$ such that  $p_{\tau(i)}=p_{\tau(j)}=p^{(k)}$ if and only if $i,j\in \mathcal A_k$, for some $1\le k\le K$. Clearly, $d_k=card\left(\mathcal A_k\right)$. Then,
\begin{equation}\label{schur}
s_{\lambda^0}(p)=\frac{\underset{\sigma\in {\cal S}_m}{\sum}(-1)^\sigma\prod_{k=1}^K\prod_{i\in \mathcal A_k}\left(p_{\tau(i)}^{m-\sigma(i)-m_k-d_k+\tau(i)}h_{\sigma(i)}^{m_k+d_k-\tau(i)}\right)}{\prod_{k=1}^K\left(0!1!\cdots\left(d_k-1\right)!\right)\prod_{k<l}\left(p^{(k)}-p^{(l)}\right)^{d_kd_l}},
\end{equation}
\noindent where ${\cal S}_m$ is the set of all the permutations of $\left\{1,...,m\right\}$ and where  $h_i=\lambda^0_i+m-i$  for $i=1,...,m$.

The next theorem gives, for inhomogeneous random words, both $\mathbb P_{m,n}(\lambda_0)$ and the distribution of $LI_\al$. The first statement is due to Its, Tracy and Widom ([\ref{ItTrWi}], [\ref{ItTrWi1}]), while the second follows directly from the fact that the length of the longest increasing subsequence is
equal to the length of the first row of the corresponding Young diagrams.

\begin{theo}\label{PushForward}
\begin{itemize}
\item[(i)] On $[m]^n$, the image (or push-forward) of $\mathbb P_{W,m,n}$ by the mapping $S: [m]^n\rightarrow \mathcal P_m^{(n)}$ is, for any $\lambda^0=(\lambda^0_1,\lambda^0_2,...,\lambda^0_m)\in\mathcal P_m^{(n)}$, given by
\begin{equation}\label{NonUniPlan}
\mathbb P_{m,n}(\lambda^0)=s_{\lambda^0}(p)f^{\lambda^0}.
\end{equation}
\noindent Above, $f^{\lambda^0}$ is the number of Young diagrams of shape $\lambda^0$ with elements in $\{1,...,n\}$:
$$f^{\lambda^0}=n!\underset{1\le i<j\le m}{\prod}(\lambda^0_i-\lambda^0_j+j-i)\prod_{j=1}^{m}\frac{1}{(\lambda^0_j+m-j)!},$$
\noindent and $s_{\lambda^0}(p)$ is the Schur function of shape $\lambda^0$ in the variable $p=\big(p_{\tau(1)},...,$ $p_{\tau(m)}\big)$ given in (\ref{schur}), with $\tau$ a permutation of $\left\{1,...,m\right\}$ corresponding to a non-increasing ordering of  $p_1,p_2,...,p_m$.%
\item[(ii)] The Poissonization of $\mathbb P_{m,n}$ is the generalized Charlier ensemble $\mathbb P_{Ch,m}^\al$ defined in (\ref{GenCharlier}). In particular, for the Poissonized word problem,
\begin{equation}\label{NonUniPoissonized}
\mathbb P_{W,m}^\al\left(LI_\al\le t\right):=e^{-\al}\sum_{n=0}^\infty\mathbb P_{m,n}\left(\lambda_1\le t\right)\frac{\al^n}{n!}=\mathbb P_{Ch,m}^\al\left(\lambda_1\le t\right).
\end{equation}
\end{itemize}
\end{theo}

For uniform alphabet, Johansson [\ref{Joha2}] obtained the convergence, as $\al\rightarrow \infty$, of the Poissonized measure on $\mathcal P_m$ to the joint law of the ordered eigenvalues of the GUE.
Next, following his lead and techniques, we generalize this result to the non-uniform case, where
the convergence is towards the joint law of the eigenvalues $\left(\xi_1,...,\xi_m\right)$, ordered within each block, of an element of $\gp$. The density of $\left(\xi_1,...,\xi_m\right)$ is, for any $x\in \R^m$, given by
\begin{equation}\label{DensityDirectSumGUE1}
f_{\xi_1,...,\xi_m}(x)=\frac{1}{\sqrt{2\pi}}c_{m}\prod_{k=1}^K\Delta_{k}(x)^2e^{-\sum_{i=1}^mx^2_i/2},
\end{equation}
\noindent where $c_{m}=(2\pi)^{-(m-1)/2}\prod_{k=1}^K\left(0!1!\cdots\left(d_k-1\right)!\right)^{-1}$, and where $$\Delta_{k}(x)=\underset{\tiny{m_k+1\le i<j\le m_k+d_k}}{\prod}\left(x_i-x_j\right).$$

\begin{theo}\label{theoPoissonized} Let $\lambda(RSK({\bf X}_W))=\left(\lambda_1,...,\lambda_m\right)$ be the common shape of the Young diagrams associated with $W$ through the RSK correspondence. Let $\left(\xi_1,...,\xi_m\right)$ be the eigenvalues of an element of $\gp$, written in such a way that $\xi_{m_k+1}\ge \cdots\ge \xi_{m_k+d_k}$ for $k=1,...,K$, and let $f_{\xi_1,...,\xi_m}$ be its density given by (\ref{DensityDirectSumGUE1}). Then, for any continuous function $g$ on $\R^m$,
\begin{equation} \underset{\al\rightarrow\infty}{\lim}\mathbb E_{m}^\al\left(g\left(\frac{\lambda_1-\al p_{\tau(1)}}{\sqrt{\al p_{\tau(1)}}},...,
\frac{\lambda_m-\al p_{\tau(m)}}{\sqrt{\al p_{\tau(m)}}}\right)\right)= \int_{\R^m}g(x)f_{\xi_1,...,\xi_m}(x)dx.
\end{equation}
\end{theo}

\noindent {\bf Proof.} By Theorem \ref{PushForward}, for any partition $\lambda^0=(\lambda^0_1,\lambda^0_2,...,\lambda^0_m)$ of $n\in \N$,
$$\mathbb P_{m,n}(\lambda(RSK({\bf X}_W))=\lambda^0)=s_{\lambda^0}(p)f^{\lambda^0},$$
\noindent where
$$f^{\lambda^0}=n!\underset{1\le i<j\le m}{\prod}(\lambda^0_i-\lambda^0_j+j-i)\prod_{j=1}^{m}\frac{1}{(\lambda^0_j+m-j)!},$$
\noindent and where $s_{\lambda^0}(p)$ is the Schur function of shape $\lambda^0$ in the variable  $p=\left(p_{\tau(1)},...,p_{\tau(m)}\right)$ as given in (\ref{schur}). Hence the Poissonized measure is
$$\mathbb P_m^\al\left(\lambda^0\right)=e^{-\al}\sum_{n=0}^\infty n!\underset{1\le i<j\le m}{\prod}(\lambda^0_i-\lambda^0_j+j-i)\prod_{j=1}^{m}\frac{1}{(\lambda^0_j+m-j)!}s_{\lambda^0}(p)\frac{\al^n}{n!}.$$
\noindent Next, for $i=1,...,m$, let $$x_i=\frac{\lambda^0_i-\al p_{\tau(i)}}{\sqrt{\al p_{\tau(i)}}},$$
\noindent then, as $\al\rightarrow \infty$,
\begin{align}
\prod_{j=1}^{m}\!\!\frac{1}{(\lambda^0_j+m-j)!}\sim (2\pi)^{-m/2}\frac{e^\al}{\al^n}\al^{-m(m-1)/2}\left(\prod_{i=1}^{m}p_{\tau(i)}^{\tau(i)-m}\right)e^{-\sum_{i=1}^mx_i^2/2},
\end{align}
\noindent and
\begin{align}
&\underset{1\le i<j\le m}{\prod}(\lambda^0_i-\lambda^0_j+j-i)\nonumber\\
&\sim \al^{m(m-1)/2-\sum_{k=1}^Kd_k(d_k-1)/4}\prod_{k=1}^K\left(\left(p^{(k)}\right)^{d_k(d_k-1)/4}\Delta_k(x)\right)
\prod_{k< l}\left(p^{(k)}-p^{(l)}\right)^{d_kd_l}.
\end{align}
\noindent Together with
\begin{align}
&\underset{\sigma\in {\cal S}_m}{\sum}(-1)^\sigma\prod_{k=1}^K\prod_{i\in \mathcal A_k}\left(p_{\tau(i)}^{m-\sigma(i)-m_k-d_k+\tau(i)}h_{\sigma(i)}^{m_k+d_k-\tau(i)}\right)\nonumber\\
&\sim \prod_{i=1}^{m}p_{\tau(i)}^{m-\tau(i)}\prod_{k=1}^{K}\left(p^{(k)}\right)^{-d_k(d_k-1)/2}\al^{\sum_{k=1}^Kd_k(d_k-1)/4}\prod_{k=1}^K\left(\left(p^{(k)}\right)^{d_k(d_k-1)/4}\Delta_k(x)\right),
\end{align}

\noindent the limiting density of $\left(\left({\lambda_1-\al p_{\tau(1)}}\right)/{\sqrt{\al p_{\tau(1)}}},...,
\left({\lambda_m-\al p_{\tau(m)}}\right)/{\sqrt{\al p_{\tau(m)}}}\right)$, as $\al\rightarrow \infty$, is
$$\sqrt{2\pi}c_{m}\prod_{k=1}^K\Delta_{k}(x)^2e^{-\sum_{i=1}^mx^2_i/2},\ \ x=(x_1,...,x_m)\in \R^m,$$
\noindent which is just the joint density of the eigenvalues, ordered within each block, of an element of $\GKM$. The statement then follows from a Riemann sums approximation argument as in [\ref{Joha2}]. \hfill $\Box$

The next result is concerned with "de-Poissonization", and again is the non-uniform version (with a similar proof) of a result of Johansson.

\begin{prop}\label{dePoissonization}
Let $\al_n=n+3\sqrt{n\log n }$ and  $\beta_n=n-3\sqrt{n\log n}$. Then there is a constant $C$ such that, for sufficiently large $n$, and for any $0\le n_i \le n$, $i=1,...,m$,
\begin{align}\label{IneDePoisson}
\mathbb P_m^{\al_n}\left(\lambda_1\le n_1,...,\lambda_m\le n_m\right)-\frac{C}{n^2}&\le\mathbb P_{m,n}\left(\lambda_1\le n_1,...,\lambda_m\le n_m\right)\nonumber\\
&\le \mathbb P_m^{\beta_n}\left(\lambda_1\le n_1,...,\lambda_m\le n_m\right)+\frac{C}{n^2}.
\end{align}

\end{prop}

\noindent {\bf Proof.} The proof is analogous to the proof of the corresponding uniform alphabet result, given in [\ref{Joha2}] (see also Lemma 4.7 in [\ref{BoOkOl}]).   First, a simple consequence of the description of
the RSK correspondence ensures that $\mathbb P_{m,n}\left(\lambda_1\le n_1,...,\lambda_m\le n_m\right)$ is
non-increasing in $n$, i.e.,
\begin{equation}\label{decrea}
\mathbb P_{m,n+1}\left(\lambda_1\le n_1,...,\lambda_m\le n_m\right) \le \mathbb P_{m,n}\left(\lambda_1\le n_1,...,\lambda_m\le n_m\right).
\end{equation}
Next,
$$\mathbb P_m^{\al}\left(\lambda_1\le n_1,...,\lambda_m\le n_m\right) = \sum_{n=0}^{\infty}e^{-\al}\frac{\al^n}{n!}\mathbb P_{m,n}\left(\lambda_1\le n_1,...,\lambda_m\le n_m\right),$$  and then, proceeding as in [\ref{Joha2}],

\begin{align}\label{esti}
\Big|\mathbb P_m^{\al}\left(\lambda_1\le n_1,...,\lambda_m\le n_m\right) - &\sum_{|n-\al|\le
\sqrt{8\al\log\al}}e^{-\al}\frac{\al^n}{n!}\mathbb P_{m,n}\left(\lambda_1\le n_1,...,\lambda_m\le n_m\right) \Big| \nonumber\\
&\le \frac{C}{\al^2},
\end{align}
for some constant $C$, $\al$ sufficiently large and all $1\le n_i \le n$, $i=1, \dots, m$.
Replacing $\al$ by respectively $n+3\sqrt{n\log n }$ and $n-3\sqrt{n\log n }$
completes the proof.\hfill$\Box$

We are now ready to obtain asymptotics for the shape of the Young diagrams associated
with a random word $W\in [m]^n$, when $m$ and $n$ go to infinity.
Before stating our result, let us recall the well known, large $m$, asymptotic behavior
of the spectrum of the $m\times m$ GUE ([\ref{r18}], [\ref{TrWi2}], [\ref{Joha2}]):

Let $\xi^{GUE,m}_j$ be the $jth$ largest eigenvalue of an element of the $m\times m$ GUE.
For each $r\ge 1$, there is a distribution function $F_r$ on $\R^r$, such that, for all $(t_1,...,t_r)\in\R^r$,
$$\lim_{m\rightarrow\infty}\mathbb P_{GUE,m}\left(\xi^{GUE,m}_j\le 2\sqrt{m}+t_j/m^{1/6}, j=1,...,r\right)=F_r(t_1,...,t_r).$$

The multivariate distribution function $F_r$ originates in [\ref{r18}] and [\ref{TrWi2}], another expression
for it is also given in [\ref{Joha2}] (see (3.48) there) and its one dimensional marginals are
Tracy-Widom distributions.

Once more, our next theorem is already present, for uniform alphabets, in Johansson [\ref{Joha2}].
\begin{theo}\label{theoTracyWidom} Let $r\ge 1$.  Let $d_1\rightarrow +\infty$, as $m\rightarrow+\infty$.
Then, for all $(t_1,...,t_r)\in\R^r$,
\begin{align}\label{ConvPoissonized}
\underset{m\rightarrow\infty}{\lim}\underset{\al\rightarrow\infty}{\lim}\mathbb P_m^\al\bigg(\lambda_j\le \al \pmax+2\sqrt{d_1\al \pmax }+t_jd_1^{-1/6}&\sqrt{\al \pmax}, j=1,...,r\bigg)\nonumber\\
&=F_r(t_1,...,t_r),
\end{align}
\noindent and,
\begin{align}\label{ConvGrowing}
\lim_{d_1\rightarrow\infty}\lim_{n\rightarrow\infty}\mathbb P_{m,n}\!\bigg(\lambda_j\le n\pmax +2\sqrt{d_1n\pmax }+t_j&d_1^{-1/6}\sqrt{n\pmax }, j=1,...,r\bigg)\nonumber\\
&=F_r(t_1,...,t_r).
\end{align}
\end{theo}

\noindent {\bf Proof.} By Theorem \ref{theoPoissonized}, for each $r\ge 1$, and for all $(s_1,...,s_r)\in \R^r$,
\begin{equation}\label{dePoissonLw}
\lim_{\al\rightarrow\infty}\mathbb P_{W,m}^\al\left(\frac{\lambda_j-\al \pmax}{\sqrt{\al \pmax}}\le s_j,\ j=1,...,r\right)= \mathbb P_{GUE,d_1}\left(\xi_j\le s_j,\ j=1,...,r\right),
\end{equation}
\noindent where $\xi_j$ is the $jth$ largest eigenvalue of the $d_1\times d_1$ GUE.  Hence, for any $(t_1,...,t_r)\in\R^r$,
\begin{align}\label{AlphaBeta}
&\lim_{\al\rightarrow\infty}\mathbb P_m^\al\left(\lambda_j\le \al \pmax+2\sqrt{d_1\al \pmax}+t_jd_1^{-1/6}\sqrt{\al \pmax}, j=1,...,r\right)\nonumber\\
&=\lim_{\al\rightarrow\infty}\mathbb P_m^\al\left(\frac{\lambda_j-\al \pmax}{\sqrt{\al \pmax}}\le2\sqrt{d_1}+t_jd_1^{-1/6}, j=1,...,r\right)\nonumber\\
&=\mathbb P\left(\xi_j\le 2\sqrt{d_1}+t_jd_1^{-1/6}, j=1,...,r\right).
\end{align}
\noindent As $d_1\rightarrow \infty$, the result of Tracy-Widom on the convergence
of the spectrum of the GUE  gives the first conclusion, proving (\ref{ConvPoissonized}).
Next, by Proposition~\ref{dePoissonization}, with $\al_n=n+3\sqrt{n\log n}$ and $\beta_n=n-3\sqrt{n\log n}$, there is a constant $C$ such that, for sufficiently large $n$, and for any $0\le s_j \le n$, $j=1,...,r$,
\begin{align}\label{dePoissonShape}
\mathbb P_m^{\al_n}\left(\lambda_j\le s_j, j=1,...,r\right)-\frac{C}{n^2}&\le\mathbb P_{m,n}\left(\lambda_j\le s_j, j=1,...,r\right)\nonumber\\
&\le \mathbb P_m^{\beta_n}\left(\lambda_j\le s_j, j=1,...,r\right)+\frac{C}{n^2}.
\end{align}
\noindent But, $n=\left(1-\varepsilon_\al\right)\al_n$, with $\varepsilon_\al=3\sqrt{n\log n}/\left(n+3\sqrt{n\log n}\right)$, whereas $n=\left(1+\varepsilon_\beta\right)\beta_n$ with $\varepsilon_\beta=3\sqrt{n\log n}/\left(n-3\sqrt{n\log n}\right)$. Since $\varepsilon_\al, \varepsilon_\beta\rightarrow 0$, as $n\rightarrow \infty$, it follows from (\ref{dePoissonShape}), by setting $s_j=n\pmax +2\sqrt{d_1n\pmax }+t_jd_1^{-1/6}\sqrt{n\pmax }$, that
\begin{align}\label{dePoissonShape2}
&\lim_{n\rightarrow\infty}\mathbb P_m^{\al_n}\left(\lambda_j\le \al_n\pmax +2\sqrt{d_1\al_n\pmax }+t_jd_1^{-1/6}\sqrt{\al_n\pmax }, j=1,...,r\right)\nonumber\\
&\le\lim_{n\rightarrow\infty}\mathbb P_{m,n}\left(\lambda_j\le  n\pmax+2\sqrt{d_1 n\pmax}+t_jd_1^{-1/6}\sqrt{ n\pmax}, j=1,...,r\right)\nonumber\\
&\le\lim_{n\rightarrow\infty} \mathbb P_m^{\beta_n}\left(\lambda_j\le \beta_n\pmax+2\sqrt{d_1\beta_n\pmax}+t_jd_1^{-1/6}\sqrt{\beta_n\pmax}, j=1,...,r\right).
\end{align}
\noindent Now, (\ref{AlphaBeta}) holds true with $\al$ replaced by $\al_n$ or $\beta_n$. Finally, (\ref{ConvGrowing}) follows from
(\ref{dePoissonShape2}) by letting $d_1\rightarrow \infty$. \hfill $\Box$

\begin{remk}  The convergence results in Theorem \ref{theoTracyWidom} are obtained by taking successive limits, i.e., first in $n$ and then in $m$. For uniform finite alphabets, in which case $d_1=m$, Johansson [\ref{Joha2}]
obtained the simultaneous convergence, for the length of the longest increasing subsequence, via a
careful analysis of corresponding kernels and methods of orthogonal polynomials. These results demand: $(\log n)^{3/2}/m\rightarrow0$ and $\sqrt{n}/m\rightarrow\infty$.
Also in the uniform case, under the assumption $m=o\left(n^{3/10}(\log n)^{-3/5}\right)$, the simultaneous convergence result (\ref{ConvGrowing}) is obtained, via Gaussian approximation, in [\ref{BrHo}] where non-uniform results are also given.
\end{remk}

\section{Appendix}

Let $\xi^{GUE,m}_{max,0}$ (resp. $\xi^{GUE,m}_{max}$) be the maximal eigenvalue of
an element of the $m\times m$ traceless GUE
(resp. GUE).  Below, we give simple proofs of the convergence of $\xi^{GUE,m}_{max,0}/\sqrt m$
(or equivalently of $\xi^{GUE,m}_{max}$) towards $2$. These proofs are based on
the "tridiagonalization" technique originating
in Trotter [\ref{Trot}] (see also Silverstein [\ref{Si}] where similar ideas are used).
Our first result is the well known Householder representation of Hermitian matrices.

\begin{lem}\label{lemma5.1}
Let ${\bf G}=(G_{i,j})_{1\le i, j\le m}$ be an element of the GUE.
Then, there exists a unitary matrix $\bf U$, such that
\begin{equation}
{\bf T}:={\bf UGU^*}=\left(\begin{matrix}
A_{1,1} & \chi_{m-1}^2 & 0 & \cdots & 0 \\
\chi_{m-1}^2 & A_{2,2} & \chi_{m-2}^2 & \cdots & 0 \\
\vdots  & \ddots & \ddots  & \ddots & \vdots \\
0  & \cdots & \chi_2^2 & A_{m-1,m-1} & \chi_1^2 \\
0 & \cdots & 0 & \chi_1^2 & A_{m,m} \\
\end{matrix}\right),
\end{equation}

\noindent where $A_{1,1},...,A_{m,m}$ are independent standard
normal random variables, and for each $1\le k\le m-1$, $\chi_{m-k}^2$ has a chi-squared distribution, with $m-k$  degrees of freedom. Moreover, for each $k=1,...,m-1$, $A_{k,k}$ is independent of $\chi_{m-k}^2,...,\chi_1^2$.
\end{lem}

\begin{prop}\label{theo4.1} Let $\xi^{GUE,m}_{max,0}$ (resp. $\xi^{GUE,m}_{max}$) be the maximal
eigenvalue of an element of the $m\times m$ traceless GUE (resp. GUE), then as ${m\rightarrow\infty}$,
$$\frac{\xi^{GUE,m}_{max,0}}{\sqrt m}\rightarrow 2, \ \ \  \left(resp.
\frac{\xi^{GUE,m}_{max}}{\sqrt m}\rightarrow 2\right) \ almost\ surely.$$
\end{prop}

\noindent {\bf Proof.}  An elementary proof is obtained along the following lines:
First, by Lemma~\ref{lemma5.1}, ${\bf G}$ and ${\bf T}$ share the same eigenvalues.
Next, by the Ger$\hat{\text{s}}$gorin circle theorem (see [\ref{r21}]), for any eigenvalue
$\xi_i$ of ${\bf G}$, letting also $\chi_0^2=\chi_m^2=0$,
$$\xi_i\in\underset{k=1,...,m}{\bigcup}\left[A_{k,k}-\chi_{m-k+1}^2-\chi_{m-k}^2,A_{k,k}+\chi_{m-k+1}^2+\chi_{m-k}^2\right].$$

\noindent Hence
\begin{align}\label{f5.88}
\frac{\xi^{GUE,m}_{max}}{\sqrt m}&\le \underset{k=1,...,m}{\max}\Bigg(\frac{A_{k,k}}{\sqrt m}+\frac{\chi_{m-k+1}^2}{\sqrt m}+\frac{\chi_{m-k}^2}{\sqrt m}\Bigg).
\end{align}

\noindent For each $k=1,...,m$, $A_{k,k}\sim N\big(0,1\big)$, and thus very classically
$\underset{k=1,...,m}{\max}{A_{k,k}}/\!{\sqrt m}$ $\stackrel{a.s.}{\rightarrow}0$.  Next, for any fixed $\varepsilon>0$,

\begin{align}
&\mathbb P\left(\left|\underset{k=1,...,m}{\max}\frac{\chi_{m-k+1}^2}{m}-1\right|>\varepsilon\right)\nonumber\\
&\le\mathbb P\left(\chi_{m}^2<m(1-\varepsilon)\right)+m \mathbb P\left(\chi_{m}^2>m(1+\varepsilon)\right),
\end{align}
and the tail behavior of $\chi_m^2$ ensures that $\sum_{m=1}^\infty m\mathbb P\left(\chi_{m}^2>m(1+\varepsilon)\right)<+\infty,$ and that $\sum_{m=2}^\infty \mathbb P\left(\chi_{m}^2<m(1-\varepsilon)\right)<+\infty$.  Therefore, $\underset{k=1,...,m}{\max}\chi_{m-k+1}^2/m\stackrel{a.s.}{\rightarrow}1$, and almost surely,
\begin{align}\label{e4.99}
\limsup_{m\rightarrow\infty}\frac{\xi^{GUE,m}_{max}}{\sqrt m}\le 2.
\end{align}
\noindent Next, since the empirical distribution of the eigenvalues $\left(\xi^{GUE,m}_{i}/\sqrt m\right)_{1\le i\le m}$ converges almost surely to the semicircle law $\nu$ with density $\sqrt{4-x^2}/2\pi$, for any $\varepsilon>0$,
\begin{equation}\label{f5.86}
\mathbb P\left(\liminf_{m\rightarrow\infty}\frac{\xi^{GUE,m}_{max}}{\sqrt m}>2-\varepsilon\right)=1.
\end{equation}
Letting $\varepsilon\rightarrow 0$ in (\ref{f5.86}) yields,
\begin{equation}\label{liminf}
\liminf_{m\rightarrow\infty}\frac{\xi^{GUE,m}_{max}}{\sqrt m}\ge 2\ \ \ a.s.
\end{equation}
\noindent Combining (\ref{e4.99}) and (\ref{liminf}), ${\xi^{GUE,m}_{max}}/{\sqrt m}\rightarrow 2$ almost surely,
and a similar result also
follows for ${\xi^{GUE,m}_{max,0}}/{\sqrt m}$. \hfill $\Box$

To prove our next convergence result, we first need a simple lemma.

\begin{lem}\label{EmaxChi2}
For each $k=1,2,...$, let $\chi^2_k$ be a chi-square random variable with $k$ degrees of freedom. Then,
\begin{equation}\label{MaxChi2}
\lim_{m\rightarrow\infty}\mathbb E\left(\frac{\underset{k=1,...,m}{\max}\chi^2_k}{m}\right)=1.
\end{equation}
\end{lem}

\noindent {\bf Proof.} First,
$$\mathbb E\left({\underset{k=1,...,m}{\max}\chi^2_k}\right)\ge \mathbb E\left({\chi^2_m}\right)=m.$$
\noindent Next, by the concavity of the logarithm, for any $0<t<1/2$,
\begin{align}\label{f4.24}
t\mathbb E\left(\frac{\underset{k=1,...,m}{\max}\chi^2_k}{m}\right)&\le {\frac{1}{m}}\ln\left(\sum_{k=1}^m\mathbb Ee^{t\chi^2_k}\right)\nonumber\\
&\le \frac{1}{m}\ln\left(m\frac{1}{(1-2t)^{m/2}}\right)\nonumber\\
&=\frac{\ln m}{m}-\frac{1}{2}\ln\left(1-2t\right).
\end{align}

\noindent Hence,
$$t\limsup_{m\rightarrow\infty}\mathbb E\left(\frac{\underset{k=1,...,m}{\max}\chi^2_k}{m}\right)\le -\frac{1}{2}\ln\left(1-2t\right),$$
\noindent and letting $t\rightarrow 0$, $$\limsup_{m\rightarrow\infty}\mathbb E\left(\frac{\underset{k=1,...,m}{\max}\chi^2_k}{m}\right)\le \lim_{t\rightarrow 0}-\frac{\ln\left(1-2t\right)}{2t}=1.$$

\noindent $\bigg($Since $-\ln(1-2t)\le2t+4t^2$, for $0\le t\le 1/3$, taking $t=\sqrt{\ln m/2m}$ in (\ref{f4.24}), will give
$\mathbb E\left({\underset{k=1,...,m}{\max}\chi^2_k}/{m}\right)\le 1+2\sqrt{2\ln m/m}$, for $m>10$. $\bigg)$
\hfill $\Box$

Again, in the uniform finite alphabet case, where $p_1=\cdots=p_m=1/m$, we have $K=1$, $d_1=m$. For $k=1,...,m$, and to keep up with the notation of [\ref{HoLi}], denote by $\tilde H_m^k$
the particular version of $\hat L^k_m$, as in (\ref{hatLk}).
Let $\big(\tilde B^1(t), \tilde B^2(t), ..., \tilde B^m(t)\big)$ be the $m$-dimensional Brownian motion having covariance matrix
\begin{equation}\label{f2.1.1}
\left(\begin{array}{clrr}
1 & \rho & \cdots & \rho \\
     \rho & 1 & \cdots & \rho \\       \vdots  & \vdots & \ddots  & \vdots \\ \rho  & \rho & \cdots &
1 \\   \end{array}\right)t,
\end{equation}
\noindent with $\rho=-1/(m-1)$. Then, for $k=1,...,m$ (see also [\ref{HoLi}], [\ref{Do}]),

$$\tilde H_m^k=\sqrt{\frac{m-1}{m}}\sup\sum_{i=1}^{m}\sum_{p=1}^{k}\big(\tilde B^i(t_{i-p+1}^p)-\tilde B^i(t_{i-p}^p)\big),$$
\noindent where the $\sup$ is taken over all the subdivisions $(t_i^p)$ of $[0,1]$ as in (\ref{f5.7555}). As a corollary to Theorem~\ref{theoNonUniform} (see also [\ref{HoLi}]), for each $m\ge 2$,
\begin{equation}\label{e4.1}
\left(\tilde H_m^1,\tilde H_m^2-\tilde H_m^1,...,\tilde H_m^m-\tilde H_m^{m-1}\right)\!\stackrel{d}{=}\!\left(\xi^{GUE,m}_{1,0},\xi^{GUE,m}_{2,0},...,\xi^{GUE,m}_{m,0}\right).
\end{equation}
\noindent Moreover, convergence in $L^1$ also holds.

\begin{prop}\label{theo5.13}  As ${m\rightarrow\infty}$,
$$\frac{\xi^{GUE,m}_{max,0}}{\sqrt m}\rightarrow 2,\ \ \ \ in\ L^1.$$
\noindent Equivalently,
$$\frac{\xi^{GUE,m}_{max}}{\sqrt m}\rightarrow 2,\ \ \ \ in\ L^1.$$
\noindent Equivalently,
$$\frac{\tilde H_m^1}{\sqrt m}\rightarrow 2,\ \ \ \ in\ L^1.$$
\end{prop}

\noindent {\bf Proof.} Note that when $p_1=\cdots=p_m=1/m$, $\mathcal L^{p_1,...,p_m}_{(s_1,...,s_m)}$, given by (\ref{f6.03}) is the empty set when $s_1<0$. Hence $\xi^{GUE,m}_{max,0}$ is nonnegative (this is actually clear from the traceless requirement). By Theorem \ref{theoNonUniform}, $\tilde H_m^1$ and $\xi^{GUE,m}_{max,0}$ are equal in distribution, and so it suffices to prove that, as ${m\rightarrow\infty}$,
\begin{equation}\label{5.46}
\frac{\mathbb E\left(\xi^{GUE,m}_{max,0}\right)}{\sqrt m}\rightarrow 2.
\end{equation}
\noindent Next, by Proposition \ref{prop6.2}, $\mathbb E\left(\xi^{GUE,m}_{max,0}\right)=\mathbb E\left(\xi^{GUE,m}_{max}\right)$. Moreover, taking expectations on both sides of (\ref{f5.88}) gives:
$$\mathbb E\left({\xi^{GUE,m}_{max}}\right)\le \mathbb E\left(\underset{k=1,...,m}{\max}{A_{k,k}}\right)
+\mathbb E\left(\underset{k=1,...,m}{\max}{\chi_{m-k+1}^2}\right)+\mathbb E\left(\underset{k=1,...,m}{\max}{\chi_{m-k}^2}\right).$$
\noindent It is well known that,
$$\mathbb E\left(\underset{k=1,...,m}{\max}A_{k,k}\right)\le \sqrt{2\ln m},$$
\noindent  while, by Lemma \ref{EmaxChi2},
$$\underset{m\rightarrow\infty}{\limsup}\ \mathbb E\left(\underset{k=1,...,m}{\max}\frac{\chi_{k}^2}{\sqrt m}\right)=1,$$

\noindent leading to
$$\underset{m\rightarrow\infty}{\limsup}\ \mathbb E\left(\frac{\xi^{GUE,m}_{max,0}}{\sqrt m}\right)\le 2.$$

\noindent Now, $\xi^{GUE,m}_{max,0}$ is nonnegative and by Proposition~\ref{theo4.1},
${\xi^{GUE,m}_{max,0}}/{\sqrt m}\rightarrow 2$, almost surely. Thus, by Fatou's Lemma,
$$\underset{m\rightarrow\infty}{\liminf}\ \mathbb E\left(\frac{\xi^{GUE,m}_{max,0}}{\sqrt m}\right)\ge \mathbb E\left(\underset{m\rightarrow\infty}{\liminf}\ \frac{\xi^{GUE,m}_{max,0}}{\sqrt m}\right)=2,$$
\noindent and so, $\lim_{m\rightarrow\infty}\mathbb E\left({\xi^{GUE,m}_{max,0}}/{\sqrt m}\right)=2.$ Using once more the fact that $\xi^{GUE,m}_{max,0}$ is nonnegative, we conclude that
$\lim_{m\rightarrow\infty}\mathbb E\left|{\xi^{GUE,m}_{max,0}}/{\sqrt m}-2\right|=0,$
\noindent and by the weak law of large number,
$\lim_{m\rightarrow\infty}\mathbb E\left|{\xi^{GUE,m}_{max}}/{\sqrt m}-2\right|=0$.\hfill $\Box$

\begin{remk}
A small and elementary tightening of the arguments of Davidson and Szarek [\ref{DaSz}] will also provide an alternative proof of Proposition~\ref{theo5.13}.
\end{remk}

\noindent {\bf Proof of Proposition~\ref{prop2.7}.} By Proposition \ref{prop6.02},
$$\underset{m_k< i\le m_k+d_k}{\max}\xi^0_i=\underset{m_k< i\le m_k+d_k}{\max}\xi_i-\sqrt {p^{(k)}}\sum_{l=1}^m\sqrt{p_{l}}{\bf X}_{l,l}.$$
Since $\underset{m_k< i\le m_k+d_k}{\max}\xi_i$ is the maximal eigenvalue of an element of the $d_k\times d_k$ GUE, with probability one or in the mean,
$\lim_{d_k\rightarrow \infty}{\underset{m_k< i\le m_k+d_k}{\max}\xi_i}/{\sqrt{d_k}}=2.$
\noindent Moreover, $\sum_{l=1}^m\sqrt{p_{l}}{\bf X}_{l,l}$ is a centered Gaussian random variable with variance
$Var\left(\sum_{l=1}^m\sqrt{p_{l}}{\bf X}_{l,l}\right)=\sum_{l=1}^mp_{l}=1.$
\noindent Hence, with probability one or in the mean, $\lim_{d_k\rightarrow \infty}{\sqrt {p^{(k)}}\sum_{l=1}^m\sqrt{p_{l}}{\bf X}_{l,l}}/{\sqrt{d_k}}=0$. \hfill$\Box$

\vspace{.3in}
\noindent {\bf Acknowledgments:} Many thanks to a referee for a careful reading of the paper.


\begin{thebibliography}{30}

\bibitem{AGZ}\label{AGZ} G.W. Anderson, A. Guionnet, O. Zeitouni, An introduction to Random Matrices. Cambridge University Press, (2009).
\bibitem{BaDeJo}\label{r41} J. Baik, P.A. Deift, K. Johansson, On the distribution of the length of the longest increasing subsequence in a random permutation. {\it J. Amer. Math. Soc.}, {\bf 12}, (1999), 1119-1178.
\bibitem{Bary}\label{Bary} Y. Baryshnikov, GUEs and queues. {\it Probab. Theor. Rel. Fields.}, {\bf 119}, (2001), 256-274.
\bibitem{Bi}\label{Bi} P. Billingsley, Convergence of probability measures, 2nd ed.. John Wiley and Sons, Inc., (1999).
\bibitem{BoOkOl}\label{BoOkOl}  A. Borodin, A. Okounkov, G. Olshanki. Asymptotics of Plancherel measures for symmetric groups. {\it J. Amer. Math. Soc.} {\bf 13} (2000), 481-515.
\bibitem{BrHo}\label{BrHo} J.-C. Breton, C. Houdr\'e,  Asymptotics for random Young diagrams
when the word length and the alphabet size simultaneously grow to infinity. {\it Bernoulli},
{\bf 16}, (2010), 471-492.
\bibitem{CG}\label{CG} G.P. Chistyakov, F. G\"otze.  Distribution of the shape of Markovian
random words.  {\it Probab. Theory Related Fields}, {\bf 129} (2004), 18-36.
\bibitem{DaSz}\label{DaSz} K. Davidson, S. Szarek. Local operator theory, random matrices and Banach spaces. {\it Handbook of the Geometry of Banach Spaces}, vol. I, North Holland, (2001), 317-366.
\bibitem{Do}\label{Do} Y. Doumerc, A note on representations of classical Gaussian matrices. {\it S\'emimaire de Probabilit\'es XXXVII.}, {\it Lecture Notes in Math., No. 1832}, Springer, Berlin, (2003),  370-384.
\bibitem{Fu}\label{Fu}  W. Fulton, Young tableaux: with applications to representation theory and geometry. Cambridge University Press, (1997).
\bibitem{GlWh}\label{GlWh} P.W. Glynn, W, Whitt, Departures from many queues in series. {\it Ann. Appl.\ Probab.} {\bf 1} (1991), 546-572.
\bibitem{GrTrWi}\label{GrTrWi} J. Gravner, J. Tracy, H. Widom, Limit theorems for height fluctuations in a class of discrete space and time growth models. {\it J. Stat. Phys.}, {\bf 102}, {Nos. 5-6}, (2001), 1085-1132.
\bibitem{Ha}\label{Ha} J.M. Hammersley, A few seedings of research. {\it Proc. Sixth Berkeley Symp. Math. Statist. and probability, vol 1}, University of california Press, (1972), 345-394.
\bibitem{HoJo}\label{r21}  R.A. Horn, C. Johnson, Topics in matrix analysis. Cambridge
University Press, (1991).
\bibitem{HoLi}\label{HoLi} C. Houdr\'e, T. Litherland, On the longest increasing subsequence
for finite and countable alphabets, in High Dimensional Probability V:
The Luminy Volume (Beachwood, Ohio, USA: Institute of Mathematical Statistics), (2009),
185-212.
\bibitem{HoLi2}\label{HoLi2} C. Houdr\'e, T. Litherland, Asymptotics for the length of the longest 
increasing subsequence of a binary Markov random word. ArXiv \# math.Pr/1110.1324, (2011). To appear in:  
Malliavin Calculus and Stochastic Analysis:  A Festschrift in honor of David Nualart.  
\bibitem{HoLi3}\label{HoLi3} C. Houdr\'e, T. Litherland, On the limiting shape of Young diagrams
associated with Markov random words. ArXiv \# math.Pr/1110.4570, (2011).
\bibitem{ItTrWi}\label{ItTrWi} A.R. Its, C. Tracy, H. Widom, Random words, Toeplitz determinants, and integrable systems. I. Random matrx models and their applications, {\it Math. Sci. Res. Inst. Publ.}, {\bf 40} Cambridge Univ. Press, Cambridge, (2001), 245-258.
\bibitem{ItTrWi1}\label{ItTrWi1} A.R. Its, C. Tracy, H. Widom, Random words, Toeplitz determinants, and integrable systems. II. Advances in nonlinear mathematics and science, {\it Phys. D.}, vol. 152-153 (2001), 199-224.
\bibitem{Joha2}\label{Joha2} K. Johansson, Discrete orthogonal polynomial ensembles and the Plancherel measure. {\it Ann. Math.} {\bf 153} (2001), 199-224.
\bibitem{Meh}\label{Meh}  M.L. Mehta, Random matrices, 2nd ed.
Academic Press, San Diego, (1991).
\bibitem{Ok}\label{Ok} A. Okounkov, Random matrices and random permutations. {\it Int. Math. Res. Not.} {\bf 2000}, no. 20, (2000), 1043-1095.
\bibitem{Silv}\label{Si} J. Silverstein, The smallest eigenvalue of a large-dimensional Wishart matrix. {\it Ann. Probab.} {\bf 13}, no. 4, (1985), 1364-1368.
\bibitem{St}\label{St}   R.P. Stanley, Enumerative Combinatorics. {\bf 2}, Cambridge
University Press, (2001).
\bibitem{TrWi}\label{r18}  C. Tracy, H. Widom, Level-spacing distribution and the
Airy kernel. {\it Comm. Math. Phys.} {\bf 159} (1994), 151-174.
\bibitem{TrWi2}\label{TrWi2} C. Tracy, H. Widom, Correlation functions, cluster functions, and spacing distributions for random matrices. {\it J. Statist. Phys.} {\bf 92}, no. 5-6, (1998), 809-835.
\bibitem{TrWi3}\label{TrWi3}  C. Tracy, H. Widom, On the distribution of the lengths of the longest inceasing monotone subsequences in random words. {\it Probab. Theor. Rel. Fields.} {\bf 119} (2001), 350-380.
\bibitem{Trot}\label{Trot}  H.F. Trotter, Eigenvalue distribution of large Hermitian matrices;
Wigner's semi-circle law and a theorem of Kac, Murdock, and Szeg{$\ddot{\text{o}}$}. {\it Adv. Math.} {\bf 54} (1984), 67-82.
\end{thebibliography}
\end{document}